\documentclass[12pt]{amsart}
\usepackage{amsfonts}

\usepackage{amsmath}
\usepackage{graphicx}

\newtheorem{thm}{Theorem}[section]
\newtheorem{thm*}{Theorem}
\newtheorem{cor}[thm]{Corollary}
\newtheorem{lem}[thm]{Lemma}

\theoremstyle{definition}

\theoremstyle{remark}

\numberwithin{equation}{section}

\def\K{{\mathbb K}}
\def\S{{\mathbb S}}
\def\R{{\mathbb R}}
\def\Z{{\mathbb Z}}

\def\grad{\nabla}

\def\lf{\left}
\def\ri{\right}

\allowdisplaybreaks

  \email{leung@math.umn.edu}
  \email{tomwan@math.cuhk.edu.hk}
\thanks{${}^*$ Partially supported by   NSF/DMS-0103355 and IMS,
CUHK.}
\thanks{${}^{**}$ Partially supported by 
   Earmarked Grants of Hong Kong CUHK4291/00P}

\input{tcilatex}

\begin{document}
\title{Harmonic maps and the topology of conformally compact Einstein manifolds}
\author{Naichung C. Leung${}^{\ast}$ and Tom Y. H. Wan${}^{\ast\ast}$}
\address{School of Mathematics, University of Minnesota, Minneapolis, MN 55455, USA
}
\address{Department of Mathematics, The Chinese University of Hong Kong, Shatin,
N.T., Hong Kong}

\begin{abstract}
We study the topology of a complete asymptotically hyperbolic Einstein
manifold such that its conformal boundary has positive Yamabe invariant. We
proved that all maps from such manifold into any nonpositively curved
manifold are homotopically trivial. Our proof is based on a Bochner type
argument on harmonic maps.
\end{abstract}

\maketitle

\markboth{Naichung Conan Leung and Tom Y. H. Wan}
{Harmonic maps and Einstein manifolds}

\setlength{\baselineskip}{20pt}

\section*{Introduction}

While studying the AdS/CFT correspondence in physics, Witten and Yau \cite
{Witten-Yau} proved that for any conformally compact Einstein manifold $%
(M^{n+1},g)$ with $\mbox{Ric}=-ng$ and a conformal infinity of positive
scalar curvature, one has $H_{n}\left( M,\mathbb{Z}\right) =0$, in
particular the conformal boundary of such manifold is connected. This result
has opened up a new and exciting direction in geometric analysis. Since
then, there are several different proofs of Witten-Yau theorem by Anderson 
\cite{Anderson}, Cai and Galloway \cite{Cai-Galloway} and X. Wang \cite
{X.Wang, X.Wang1}.

In terms of homotopy theory, the above result can be expressed as 
\begin{equation*}
\left[ (M,\partial M),({\S}^{1},\ast)\right] =1, 
\end{equation*}
i.e., any continuous map from $M$ to ${\S}^{1}$ which maps the boundary into
the marked point is homotopic to the constant map. So it is natural, from the point of view of harmonic maps, to ask
whether the same is true if we replace ${\S}^{1}$ by any nonpositively curved
manifold. And we answered this question affirmatively. More
precisely, we have

\begin{thm*}
\label{thm-top} Let $(M^{n+1},g)$, $n\geq 2$, be a conformally compact
Einstein manifold of order $C^{3,\alpha}$ with $\mbox{Ric}_{g}= -ng$ such
that the conformal infinity of $M$ has positive Yamabe invariant. Suppose
that $N$ is a compact nonpositively curved manifold. Then the homotopy
classes $[(M,\partial M),(N,\ast)]$ are trivial.
\end{thm*}

In case that $\mbox{Ric}\geq 0$, the usual Bochner technique on harmonic
maps is a standard tool in proving similar type of theorems. When $\mbox{Ric}%
=-ng$ there is a Matsushima formula, a version of Bochner type formula,
which can be used to obtain vanishing results under special situations. For
example, Jost-Yau \cite{Jost-Yau} and Mok-Siu-Yeung \cite{Mok-Siu-Yeung} 
had used this idea to obtain superrigidity results. In \cite{X.Wang, X.Wang1}, 
X. Wang showed that the lower bound on the smallest $L^{2}$%
-eigenvalue can be used to balance the negativity of the Ricci curvature in
the Bochner argument and hence a new proof of Witten-Yau's theorem via the
result on the smallest eigenvalue by J. Lee \cite{Lee}. Moreover, he found
the sharp lower bound on the eigenvalue such that the method works. We want
to remark that Lee's result is an estimate on the eigenvalue by the boundary
geometry. Hence, this method is in fact using the boundary geometry to
absorb the Ricci curvature term. This idea is also first introduced and used
by Witten and Yau in \cite{Witten-Yau}. Recently, P. Li and J. Wang \cite{Li-Wang}
are able to show that without assuming the manifold is conformally
compact, the conditions $\mbox{Ric}\geq -n$ and $\lambda \ge n-1$ are
sufficient to prove that the manifold has only one end of infinite volume or it is a wrap product.

We extend X. Wang's arguments to harmonic maps and generalize his result in
this paper. Suppose that $\left( M^{n+1},g\right) $ is a conformally compact
manifold of order $C^{k,\alpha}$, i.e. there exists a smooth defining
function $t$ for $\partial M$ on $\bar{M}$ such that $\bar{g}%
_{ij}=t^{2}g_{ij}$ defines a metric on $\bar{M}$ which is $C^{k,\alpha}$ up
to the boundary. By a defining function $t$, we mean $t>0$ in $M$ and $t$
vanishes to first order on $\partial M$. If $g$ satisfies the Einstein
equation $\mbox{Ric}=-ng$ and regular enough, then $|\nabla t|_{\bar{g}}=1$
on $\partial M$. Since the sectional curvature is asymptotic to $-|\nabla
t|_{\bar{g}}^{2}$ on $\partial M$ for conformally compact metric \cite
{Mazzeo}, the Einstein metric $g$ is asymptotically hyperbolic according to
Lee's terminology (called weakly asymptotically hyperbolic by X. Wang). By
denoting the first $L^{2}$-eigenvalue of $M$ as $\lambda_{g}$, we can now
state our main results in which $M$ is not necessary Einstein.

\begin{thm*}
\label{thm-main} Suppose that $(M^{n+1},g)$, $n\geq 2$, is an asymptotically
hyperbolic conformally compact manifold of order $C^{1}$ such that 
\begin{equation*}
\mbox{Ric}_{g}\geq-ng\text{ and }\lambda_{g}\ge n-1. 
\end{equation*}
Suppose that $f:M\rightarrow N$ is a smooth harmonic map of finite total energy from $M$ into a
complete nonpositively curved manifold $N$. If $\lambda_g > n-1$, then $f$ is a constant map. If $\lambda_g=n-1$, then either $f$ is a constant map or $M=\R\times \Sigma$ with the warped product metric $g=dt^2+\cosh ^2 (t) h$, where $(\Sigma, h)$ is a compact manifold with $\mbox{Ric}_h\ge -(n-1)$.
\end{thm*}
After finished the first draft of this paper, the second author was told by P. Li that using the method in \cite{Li-Wang}, theorem \ref{thm-main} remains true without assuming the conformally compactness of the manifold $M$ if the harmonic map is assumed to be asymptotically constant at infinity.

Combining theorem \ref{thm-main} with results of Li-Tam \cite{Li-Tam} and Liao-Tam 
\cite{Liao-Tam} concerning harmonic map heat flow, we can prove the
corresponding topological result for such a manifold easily. Then finally,
together with Lee's result \cite{Lee}, we can generalize Witten-Yau theorem
as stated precisely in theorem \ref{thm-top}.

To prove theorem \ref{thm-main}, we proceed as in \cite{X.Wang, X.Wang1} and consider 
the $(n-1)/2n$-power of the
energy density. The technical part is to show that under our assumption,
this power of energy density still belongs to $L^{1,2}(M)$. Our method is to show that the energy decay sufficiently fast as precisely state in lemma \ref{lem-1}. This decay lemma may be interesting by itself.

This paper is arranged as follow. In section \ref{sec-decay}, we first study
the decay rate of the energy for harmonic maps with finite total energy on
each end and prove couple of lemmas which are needed in the proof of our
theorems. We will prove our main theorem in section \ref{sec-van} and apply
it in section \ref{sec-app} to prove the desired and other topological
theorems.

Acknowledgment:  The first author would like to thank J.P. Wang for useful
discussions. The second author would like to thank P. Li, R. Schoen and L-F. Tam
for helpful discussions, X. Wang for sending us his papers and thesis, and
also T. Au, K. S. Chou, J. Wolfson and D. Pollack for their interest
in the problem. This joint work was started during the first author's visit to the Institute
of Mathematical Science of the Chinese University of Hong Kong and we would like
to thank the institute for giving us this chance to work together.

\section{Energy decay of harmonic maps with finite total energy}

\label{sec-decay}

Let $(M^{n+1},g)$ be an asymptotically hyperbolic conformally compact
manifold of order $C^{1}$. Its boundary components and the corresponding
ends will be denoted by $\Sigma_{0}^{i}=\partial M^{i}$ and $E^{i}$
respectively. Under our assumption, there exists \emph{special} defining
function $t$ in the following sense \cite{Graham-Lee}: each end $E^{i}$ can
be parametrized by $(t,x)$ such that the metric can be written as 
\begin{equation*}
g=t^{-2}\left( dt^{2}+h(t,x)\right) , 
\end{equation*}
(i.e. $\bar{g}=dt^{2}+h$) where $x\in\Sigma_{0}^{i}$, $t>0$ is small and $%
h(t,\cdot)$ is a family of metric defined on $\Sigma_{0}^{i}$.

For any $0<t_{1}<t_{2}$ sufficiently small, we define 
\begin{equation*}
E_{t_{1}}^{i}=\{p\in E^{i}:0<t(p)<t_{1}\}, 
\end{equation*}
\begin{equation*}
A_{t_{1},t_{2}}^{i}=\{p\in E^{i}:t_{1}<t(p)<t_{2}\}, 
\end{equation*}
\begin{equation*}
\Sigma_{t_{1}}^{i}=\{p\in E^{i}:t(p)=t_{1}\}, 
\end{equation*}
and 
\begin{equation*}
M_{t_{1}}=M\setminus\cup_{i}E_{t_{1}}^{i}. 
\end{equation*}
In order to simplify notations, we will omit the superscript $i$ and write $%
E_{t_{1}}$, $A_{t_{1},t_{2}}$, or $\Sigma_{t_{1}}$ if we are working on a
fixed end.

\begin{lem}
\label{lem-1} Suppose that $M^{n+1}$, $n\ge 2$, is a conformally compact manifold of
order $C^{1}$ with $\mbox{Ric}\geq-n$ and $N$ is a complete nonpositively
curved manifold. Let $f:E_{t_{0}}\rightarrow N$ be a smooth harmonic map
from an end $E_{t_{0}}\subset M$ into $N$ of finite total energy. Then 
\begin{equation*}
\int_{E_{t}}|\nabla f|^{2}=O\left( t^{n}\right) \quad\mbox
{as }t\rightarrow0. 
\end{equation*}
\end{lem}

\noindent\emph{Proof.}
We may assume that $t_{0}>0$ is sufficiently small so that $%
E_{t_{0}}$ is parametrized by $(t,x)$, where $x$ belongs to the
corresponding boundary component $\Sigma_{0}$ and $t\in(0, t_{0} )$. 

By straight forward calculations and the harmonicity of $f$, one has the following conservation law: for any smooth vector field $X$ on $E_{t_{0}}$,
\begin{equation*}
\int_{\partial A_{\tau,t}}\frac{1}{2}|\grad f|^2 \langle X,n\rangle d\sigma= 
\int_{\partial A_{\tau,t}} \langle df(X),df(n)\rangle d\sigma + 
\int_{A_{\tau,t}} \langle S_f,\grad X\rangle,
\end{equation*}
where $S_f=\frac{1}{2}|\grad f|^2g-f^*ds_N^2$ is the stress-energy tensor of $f$. There is a distinguish vector field $X=-t\frac{\partial}{\partial t}$ on $E_{t_{0}}$ and the conservation law applied to this $X$ gives
\begin{eqnarray*}
\lefteqn{ \int_{\Sigma_{\tau}}\frac{1}{2}|\grad f|^2 \lf\langle -\tau\frac{\partial}{\partial t},-\tau\frac{\partial}{\partial t}\ri\rangle d\sigma_{\tau} +
\int_{\Sigma_{t}}\frac{1}{2}|\grad f|^2 \lf\langle -t\frac{\partial}{\partial t},t\frac{\partial}{\partial t}\ri\rangle d\sigma_{t}}\\
&=& \int_{A_{\tau,t}}\langle S_f,\grad X\rangle + \int_{\Sigma_{\tau}}\lf\langle df\lf( -\tau\frac{\partial}{\partial t}\ri),df\lf( -\tau\frac{\partial}{\partial t}\ri)\ri\rangle d\sigma_{\tau} \\
& & \mbox{}+ \int_{\Sigma_{t}}\lf\langle df\lf(-t\frac{\partial}{\partial t}\ri),df\lf( t\frac{\partial}{\partial t}\ri)\ri\rangle d\sigma_{t} .
\end{eqnarray*}
That is
\begin{eqnarray}\label{eqn-con1}
\lefteqn{ \frac{1}{2}\int_{\Sigma_{\tau}} |\grad f|^2 d\sigma_{\tau} -
\frac{1}{2}\int_{\Sigma_{t}} |\grad f|^2 d\sigma_{t}} \nonumber \\
& =&
\int_{\Sigma_{\tau}}\tau^2\lf| \frac{\partial f}{\partial t} \ri|^2 d\sigma_{\tau} -
\int_{\Sigma_{t}}t^2\lf| \frac{\partial f}{\partial t} \ri|^2 d\sigma_{t}
+\int_{A_{\tau,t}}\langle S_f,\grad X\rangle.
\end{eqnarray}
For a fixed point $p=(t,x)\in E_{t_0}$, one can choose a normal coordinates $x^i$, $i=1,\ldots, n$, of the metric $h(t,\cdot)$ centered at $p$. Then $e_0=t\frac{\partial}{\partial t}$, $e_i=t\frac{\partial}{\partial x^i}$, $i=1,\ldots,n$, forms an orthonormal basis at $p$. By straight forward calculations,
$$
\grad_{e_0}X \equiv 0 \quad\mbox{and}\quad \grad_{e_i}X=e_i-\frac{t}{2}h^{jk}\frac{\partial h_{ki}}{\partial t} e_j.
$$
Therefore, at $p$,
$$
\mbox{div}X=n-\frac{t}{2}\mbox{Tr}_h \lf( h^{-1}\frac{\partial h}{\partial t} \ri)
$$
and
$$
\lf\langle \grad_{e_i}X, e_j\ri\rangle =\delta_{ij} -\frac{t}{2}\frac{\partial h_{ij}}{\partial t}
$$
for $i,j=1,\ldots,n$.
Hence at $p$,
\begin{eqnarray}\label{eqn-stress}
\lf\langle S_f,\grad X\ri\rangle &=& \sum_{A,B=0}^n S_f(e_A,e_B)\lf\langle \grad_{e_A}X, e_B\ri\rangle \nonumber \\
&=& \frac{1}{2}|\grad f|^2 \mbox{div}X-\sum_{A,B=0}^n \langle df(e_A),df(e_B)\rangle \lf\langle \grad_{e_A}X, e_B\ri\rangle \nonumber \\
&=& \frac{1}{2}|\grad f|^2 \lf[ n-\frac{t}{2}\mbox{Tr}_h \lf( h^{-1}\frac{\partial h}{\partial t} \ri) \ri] \nonumber \\
& &\mbox{} -\sum_{i,j=0}^n \langle df(e_i),df(e_j)\rangle \lf( \delta_{ij} -\frac{t}{2}\frac{\partial h_{ij}}{\partial t} \ri) \nonumber \\
&=& \frac{1}{2}|\grad f|^2 \lf[ n-\frac{t}{2}\mbox{Tr}_h \lf( h^{-1}\frac{\partial h}{\partial t} \ri) \ri] \nonumber \\
& &\mbox{}-\sum_{i=1}^n |df(e_i)|^2 + \frac{t}{2} \langle df(e_i),df(e_j)\rangle \frac{\partial h_{ij}}{\partial t}.
\end{eqnarray}
By the assumption on $g$, we first observe that
$$
|\langle S_f,\grad X\rangle | \le C|\grad f|^2
$$
and hence $|\langle S_f,\grad X\rangle |$ is integrable since $f$ has finite total energy. So we can let $\tau$ in equation (\ref{eqn-con1}) tending to 0. As in \cite{X.Wang,X.Wang1}, we can also  choose a sequence of $\tau$'s tending to $0$ such that 
$$
\int_{\Sigma_{\tau}}|\grad f|^2 \to 0 \quad \mbox{as } \tau\to 0.
$$
So, we conclude that for all $t\in (0,t_0)$,
\begin{equation}\label{eqn-con2}
\int_{E_t}\langle S_f,\grad X\rangle = \int_{\Sigma_t}t^2\lf| \frac{\partial f}{\partial t} \ri|^2 d\sigma_t -\frac{1}{2}\int_{\Sigma_t}|\grad f|^2d\sigma_t.
\end{equation}
Secondly, the assumption on $g$ and (\ref{eqn-stress}) imply that there exist $t_1\in (0,t_0)$ and $C_1>0$ such that for any $t\in (0, t_1)$,
\begin{eqnarray*}
\langle S_f,\grad X\rangle &\ge & \frac{n-C_1 t}{2}|\grad f|^2 - (1+C_1 t) \sum_{i=1}^n |df(e_i)|^2 \\
&=& \frac{n-C_1 t}{2} t^2 \lf| \frac{\partial f}{\partial t} \ri|^2 + \frac{n-2-3C_1 t}{2}\sum_{i=1}^n |df(e_i)|^2.
\end{eqnarray*}
To simplify notation, we write
$$
|\grad_t f|^2=t^2 \lf| \frac{\partial f}{\partial t} \ri|^2 \quad\mbox{and}\quad
|\grad_x f|^2=\sum_{i=1}^n |df(e_i)|^2.
$$
Putting it into (\ref{eqn-con2}), we have
\begin{eqnarray}\label{eqn-de1}
\frac{n-C_1 t}{2}\int_{E_t}|\grad_t f|^2 &+& \frac{n-2-3C_1 t}{2}\int_{E_t} |\grad_x f|^2 \nonumber \\
&\le& \frac{1}{2}\int_{\Sigma_t}t^2|\grad_tf|^2 d\sigma_t -\frac{1}{2}\int_{\Sigma_t}|\grad_x f|^2d\sigma_t.
\end{eqnarray}
If we further write
$$
F(t)=\int_{E_t}|\grad_t f|^2\quad\mbox{and}\quad G(t)=\int_{E_t} |\grad_x f|^2
$$
and note that
$$
tF'(t)=\int_{\Sigma_t}|\grad_t f|^2d\sigma_t\quad\mbox{and}\quad tG'(t)=\int_{\Sigma_t} |\grad_x f|^2d\sigma_t,
$$
then the above inequality becomes
$$
(n-C_1 t)F+(n-2-3C_1 t)G \le tF'-tG'.
$$
This immediately implies
$$
\frac{d}{dt}\lf[ t^{n-2}e^{-3C_1 t}\lf( F-G  \ri)\ri] \ge
(2n-2-4C_1 t)t^{n-3}e^{-3C_1 t}F.
$$
So if we set $t_2=\min \{ t_1, 1/2C_1 \}$, then for all $t\in (0,t_2)$,
$$
\frac{d}{dt}\lf[ t^{n-2}e^{-3C_1 t}\lf( F-G  \ri)\ri] \ge 0.
$$
As $f$ has finite total energy, $\lim_{t\to 0}F(t)=\lim_{t\to 0}G(t)=0$ and hence
$$
F(t)\le G(t)\quad \forall \, t\in (0,t_2).
$$
Put this back into (\ref{eqn-de1}), we first have
$$
tF'\ge (n-4C_1 t) F \quad \mbox{for }n=2.
$$
And then for $n\ge 3$, one simply drops the term that involve $G$ and conclude that
$$
tF'\ge (n-C_1 t) F\ge (n-4C_1 t) F.
$$
Hence in all cases,
$$
\frac{d}{dt}\lf( t^{-n}e^{4C_1 t} F(t) \ri) \ge 0.
$$
Integrating this inequality from $t$ to $t_2$, we see that there is a constant $C>0$ such that
$$
F(t)\le Ct^n \quad \forall\, t\in (0,t_2).
$$
Together with $G\le F$, we have
$$
\int_{E_t}|\grad f|^2 \le 2C t^n
$$
and the proof of the lemma is completed.

\begin{lem}
\label{lem-2} Let $f$ be a smooth harmonic map from an $(n+1)$-dimensional
manifold, then 
\begin{equation*}
|\nabla^{2}f|^{2}\geq\left( 1+\frac{1}{n}\right) |\nabla|\nabla f||^{2}. 
\end{equation*}
\end{lem}

\noindent\emph{Proof.} With respect to normal coordinates, we have $\sum_{i}
f^{\alpha}_{ii}=0$ at the center. Then for each $\alpha$, one can show as in
the gradient estimate in \cite{Schoen-Yau} that 
\begin{equation*}
|\nabla^{2} f^{\alpha}|^{2} \ge\left( 1+\frac{1}{n}\right) |\nabla|\nabla
f^{\alpha}||^{2}. 
\end{equation*}
Therefore, 
\begin{align*}
\left| \nabla|\nabla f|\right| ^{2} & = \left| \nabla\sqrt{\sum_{\alpha}
|\nabla f^{\alpha}|^{2} }\right| ^{2} \\
& = \left| \frac{\sum_{\alpha} |\nabla f^{\alpha}| \nabla|\nabla f^{\alpha}|%
}{\sqrt{\sum_{\alpha} |\nabla f^{\alpha}|^{2} }} \right| ^{2} \\
& \le \sum_{\alpha} \left| \nabla|\nabla f^{\alpha}| \right| ^{2} \\
& \le \left( 1+\frac{1}{n} \right) ^{-1} \sum_{\alpha} \left| \nabla^{2}
f^{\alpha} \right| ^{2} \\
& = \left( 1+\frac{1}{n} \right) ^{-1} \left| \nabla^{2} f \right| ^{2}.
\end{align*}

\begin{lem}
\label{lem-3} Let $M$, $N$ be as in lemma \ref{lem-1} and $f:M\rightarrow N$
be a smooth harmonic map with finite total energy, then 
\begin{equation*}
\int_{M}\left| \nabla|\nabla f|\right| ^{2}<+\infty. 
\end{equation*}
\end{lem}

\noindent\emph{Proof.} We only need to show this on one of the end $E_{t_{1}}
$. Let $\eta$ be a smooth cutoff function such that $0\leq\eta\leq1$, $%
\eta(\beta)\equiv1$ for $\beta\in(0,1]$, $\mbox{supp}\eta\subset\lbrack0,2)$%
, and $|\eta^{\prime}|\leq C$ for some absolute constant $C>0$. Then for any 
$r>0$, consider the function on $E_{t_{1}}$ defined by 
\begin{equation*}
\phi(x,t)=\eta\left( \frac{1}{r}\left| \log\frac{t_{1}}{t}\right| \right) . 
\end{equation*}
Then $0\leq\phi\leq1$, $\phi(\cdot,t)\equiv1$ for all $t\in\lbrack
e^{-r}t_{1},t_{1}]$, $\phi(\cdot,t)\equiv0$ for all $t\in(0,e^{-2r}t_{1}]$,
and $|\nabla\phi|^{2}=t^{2}\displaystyle\left( \frac{\partial\phi}{\partial t%
}\right) ^{2}\leq\displaystyle\frac{C^{2}}{r^{2}}$. Then multiple $\phi^{2}$
to the Bochner formula of $f$ and integrate, one obtains 
\begin{equation*}
\frac{1}{2}\int_{E_{t_{1}}}\phi^{2}\Delta|\nabla
f|^{2}\geq\int_{E_{t_{1}}}\phi^{2}|\nabla^{2}f|^{2}-n\int_{E_{t_{1}}}%
\phi^{2}|\nabla f|^{2}. 
\end{equation*}
By lemma \ref{lem-2}, we have 
\begin{align*}
\left( 1+\frac{1}{n}\right) \int_{E_{t_{1}}}\phi^{2}\left| \nabla|\nabla
f|\right| ^{2} & \leq n\int_{E_{t_{1}}}\phi^{2}|\nabla f|^{2}-2\int
_{E_{t_{1}}}\phi|\nabla f|\nabla\phi\cdot\nabla|\nabla f| \\
& \mbox{}+\frac{t_{1}}{2}\int_{\Sigma_{t_{1}}}\left| \frac{\partial}{%
\partial t}|\nabla f|^{2}\right| \\
& \leq n\int_{E_{t_{1}}}\phi^{2}|\nabla
f|^{2}+\int_{E_{t_{1}}}\phi^{2}\left| \nabla|\nabla f|\right| ^{2} \\
& \mbox{}+\int_{E_{t_{1}}}|\nabla\phi|^{2}|\nabla f|^{2}+\frac{t_{1}}{2}%
\int_{\Sigma_{t_{1}}}\left| \frac{\partial}{\partial t}|\nabla f|^{2}\right|
.
\end{align*}
This implies 
\begin{equation*}
\frac{1}{n}\int_{A_{e^{-r}t_{1},t_{1}}}\left| \nabla|\nabla f|\right|
^{2}\leq\left( n+\frac{C^{2}}{r^{2}}\right) \int_{E_{t_{1}}}|\nabla f|^{2}+%
\frac{t_{1}}{2}\int_{\Sigma_{t_{1}}}\left| \frac{\partial}{\partial t}|\nabla
f|^{2}\right| . 
\end{equation*}
Letting $r\rightarrow\infty$, we conclude that $\int_{E_{t_{1}}}\left|
\nabla|\nabla f|\right| ^{2}<+\infty$ which completes the proof of the lemma.

\section{Vanishing Theorem of Harmonic Maps}

\label{sec-van}

We'll prove theorem \ref{thm-main} in this section. As we mentioned in the introduction, we proceed as in \cite{X.Wang, X.Wang1} and consider the $(n-1)/2n$-power of the
energy density. The technical part is to show that under our assumption,
this power of the energy density still belongs to $L^{1,2}(M)$.

\noindent\textbf{Proof of theorem \ref{thm-main}.} Let's consider the function 
\begin{equation*}
\zeta=|\nabla f|^{\beta}. 
\end{equation*}
We first prove that for $\beta>1/2$, $\zeta\in L^{2}(M)$.

For each fixed end $E_{t_{1}}$, we set $t_{k}=2^{-k+1}t_{1}$ and consider
the integral of $\zeta^{2}$ on the annulus $A_{t_{k+1},t_{k}}$. Then H\"older
inequality implies that 
\begin{align*}
\int_{A_{t_{k+1},t_{k}}} \zeta^{2} & = \int_{A_{t_{k+1},t_{k}}} 
 |\nabla f|^{2\beta} dv_g \\
& \le \left( \int_{A_{t_{k+1},t_{k}}} 
 |\nabla f|^{2}   dv_g \right) ^{\beta} \left( \int_{A_{t_{k+1},t_{k}}} 
  dv_g \right) ^{1-\beta}  \\
& \le  \left( \int_{A_{t_{k+1},t_{k}}} |\nabla f|^{2}   dv_g \right) ^{\beta} \left( \sup_{0<t<t_1}\mbox{vol}_h(\Sigma_{t}) \int_{t_{k+1}}^{t_{k}} 
\frac{dt}{t^{n+1}} \right) ^{1-\beta} \\
& \le \lf(\frac{\sup_{0<t<t_1}\mbox{vol}_h(\Sigma_{t})}{n}\ri)^{1-\beta} \left( \int_{A_{t_{k+1},t_{k}}} |\nabla f|^{2} \right)
^{\beta} \cdot t_{k+1}^{-n(1-\beta)}.
\end{align*}
By lemma \ref{lem-1},
\begin{eqnarray*}
\int_{A_{t_{k+1},t_{k}}} \zeta^{2} &\le& C\lf(\frac{\sup_{0<t<t_1}\mbox{vol}_h(\Sigma_{t})}{n}\ri)^{1-\beta} t_{k}^{n\beta} \cdot t_{k+1}^{-n(1-\beta)} \\
& \le & C_{1} \lf(\frac{\sup_{0<t<t_1}\mbox{vol}_h(\Sigma_{t})}{n}\ri)^{1-\beta} \left( 2^{-\gamma} \right)
^{k}, 
\end{eqnarray*}
for some constants $C$ and $C_{1}>0$ independent of $k$ and $\beta$, and 
\begin{equation*}
\gamma=n\beta-n(1-\beta)=n(2\beta -1) >0. 
\end{equation*}
Hence for all $K>0$,
\begin{equation*}
\int_{A_{t_{K+1},t_{1}}} \zeta^{2} \le C_{1} \lf(\frac{\sup_{0<t<t_1}\mbox{vol}_h(\Sigma_{t})}{n}\ri)^{1-\beta}
\frac{1}{1-2^{-\gamma}}.
\end{equation*}
Since the upper bound is independent of $K$,  we've proved that
\begin{equation}\label{eqn-int}
\int_{E_{t_{1}}} \zeta^{2} \le C_{1} \lf(\frac{\sup_{0<t<t_1}\mbox{vol}_h(\Sigma_{t})}{n}\ri)^{1-\beta}
\frac{1}{1-2^{-\gamma}}.
\end{equation} 
As the end $%
E_{t_{1}}$ is arbitrary, this proves that $\zeta\in L^{2}(M)$.

Now we consider a
differential inequality of $\zeta$ which is a direct consequence of the
Bochner formula of $f$ and lemma \ref{lem-2}. To simplify notation, we use 
$\K_N\le 0$ to denote the curvature term of $N$ in the Bochner formula, which is in the order of $|\grad f|^4$. 
In this notation, the inequality becomes 
\begin{equation*}
\Delta\zeta^{2} \ge\frac{2}{\beta}\left( \frac{1}{n}+2\beta-1 \right)
|\nabla\zeta|^{2} -2n\beta\zeta^{2} - \zeta^{-\frac{2(1-\beta)}{\beta}}\K_N. 
\end{equation*}
With this inequality and $\zeta\in L^{2}(M)$, we conclude as in the proof of lemma 
\ref{lem-3} that $|\nabla\zeta|^{2}$ is also integrable. Hence $\zeta\in
L^{1,2}(M)$.

Integrating by part of the above inequality on $M_{t}$, one obtains 
\begin{eqnarray*}
-\int_{\cup_{i}\Sigma_{t}^{i}}t\frac{\partial}{\partial t}\zeta^{2}d\sigma
_{t} &\geq & \frac{2}{\beta}\left( \frac{1}{n}+2\beta-1\right)
\int_{M_{t}}|\nabla\zeta|^{2}-2n\beta\int_{M_{t}}\zeta^{2} \\
& &\mbox{ } -\int_{M_t}\zeta^{-\frac{2(1-\beta)}{\beta}}\K_N. 
\end{eqnarray*}
The left hand side can be estimated as follows: 
\begin{equation*}
\left| \int_{\cup_{i}\Sigma_{t}^{i}}t\frac{\partial}{\partial t} \zeta
^{2}d\sigma_{t}\right| \leq2
\int_{\cup_{i}\Sigma_{t}^{i}}\zeta^{2}d\sigma_{t}
+2\int_{\cup_{i}\Sigma_{t}^{i}}|\nabla\zeta|^{2}d\sigma_{t}. 
\end{equation*}
Therefore, since $\zeta\in L^{1,2}(M)$, one can choose a sequence of $t$'s
tending to $0$ as in \cite{X.Wang, X.Wang1} such that 
\begin{equation*}
\left| \int_{\cup_{i}\Sigma_{t}^{i}}t\frac{\partial}{\partial t}\zeta
^{2}d\sigma_{t}\right| \rightarrow0. 
\end{equation*}
As a consequence, 
\begin{equation}\label{eqn-flat}
-\int_{M}\zeta^{-\frac{2(1-\beta)}{\beta}}\K_N +
\frac{2}{\beta}\left( \frac{1}{n}+2\beta-1\right) \int_{M}|\nabla\zeta
|^{2}\leq2n\beta\int_{M}\zeta^{2}. 
\end{equation}
If $f$ is not a constant map, then $\zeta\not \equiv0$ and hence 
$$
\lambda_{g} \leq\frac{n\beta^{2}}{\frac{1}{n}+2\beta-1}.
$$
Since  $(n-1)/n>1/2$ for $n\ge 3$, by letting $\beta=(n-1)/n$ in this case, we have
$$
\lambda_g\le n-1.
$$
Therefore, if $\lambda_g >n-1$, then $f$ must be a constant map. If $\lambda_g=n-1$, then inequality (\ref{eqn-flat}) implies that $N$ is flat. Then the result follows from X. Wang's result \cite{X.Wang, X.Wang1}.

If $n=2$, Cheng's result \cite{Cheng} implies that $\lambda_g=n-1=1$. Hence
$$
\int_M|\grad\zeta|^2 \ge \int_M \zeta^2.
$$
Putting this into (\ref{eqn-flat}), one obtains So by taking $\beta=\frac{1}{2}+\delta$ with any $\delta>0$
$$
-\int_{M}\zeta^{-\frac{2(1-\beta)}{\beta}}\K_N \le \frac{8\delta^2}{1+2\delta} \int_M\zeta^2.
$$
To estimate $\int_M\zeta^2$, we fix a sufficiently small $t_1>0$ and note that
$$
\int_M\zeta^2 = \int_{M_{t_1}}\zeta^2+ \int_{\cup_{i}E^{(i)}_{t_1}}\zeta^2.
$$
By sub-mean value property, $|\grad f|^2$ is uniformly bounded and hence there is constant $C>0$ such that
$$
\int_{M_{t_1}}\zeta^2 =\int_{M_{t_1}}|\grad f|^{1+2\delta} \le C
$$
for sufficiently small $\delta>0$. For the second term, we apply (\ref{eqn-int}) to each $E^{(i)}_{t_1}$ and get
$$
\int_{\cup_{i}E^{(i)}_{t_1}}\zeta^2 \le C_1 V^{\frac{1}{2}-\delta}
\frac{1}{1-2^{-4\delta}},
$$
where $V=\sup_i\lf( \sup_{0<t<t_1} \mbox{vol}_h(\Sigma^{(i)}_t)/2\ri)$.
Since $1/(1-2^{-4\delta})\to +\infty$ as $\delta\to 0$, we have
$$
-\int_{M}\zeta^{-\frac{2(1-\beta)}{\beta}}\K_N \le C_2 V^{\frac{1}{2}-\delta} 
\frac{1}{1-2^{-2\delta}}\cdot\frac{8\delta^2}{1+2\delta},
$$
where $C_2$. As the upper bound tends to $0$ as $\delta\to 0$, we conclude again that $N$ is flat and hence the result follows as above. This completes the proof of the theorem \ref{thm-main}.

\section{Application to asymptotically hyperbolic Einstein manifolds}

\label{sec-app}

As an application of our vanishing theorem, we first prove the following
topological theorem.

\begin{thm}
\label{thm-top1} Suppose that $(M^{n+1},g)$, $n\geq 2$, is an asymptotically
hyperbolic conformally compact manifold of order $C^{1}$ such that 
\begin{equation*}
\mbox{Ric}_{g}\geq-ng \quad\mbox{and}\quad\lambda_{g}\ge n-1 
\end{equation*}
Then for any nonpositively curved compact manifold $N$, the homotopy classes
in $[(M,\partial M),(N,\ast)]$ are trivial or $M$ slpits as $\R\times \Sigma$ for some compact manifold $\Sigma$.
\end{thm}

\noindent\emph{Proof.} Since $M$ is conformally compact, $(M,\partial M)$ is
homotopy equivalent to $(M\setminus(\cup E_{t}^{i}),\cup\Sigma_{t}^{i}) $
for sufficiently small $t $. Therefore, in each class, we can find a smooth
representative $g:M\rightarrow N$ which maps each end $E_{t}^{i}$ into the
marked point $\ast$. Thus $g$ has bounded energy density and bounded image.
Hence by \cite{Li-Tam}, the harmonic map heat flow has a unique solution.
Also the square norm of the tension field of $g$ is in $L^{p}$ for $p>1$ and
tending to 0 near the boundary, the results in \cite{Li-Tam} imply that the
heat flow converges to a harmonic map $f$ with same boundary data as $g$. So 
$f$ is also a representative of the class $[g]$. (In particular, $f$ is
asymptotically constant on each end.) On the other hand, using the result in 
\cite{Liao-Tam}, $g$ has finite total energy (and the uniqueness) implies
that $f$ also has finite total energy. Therefore, theorem \ref{thm-main}
immediately implies that either $M$ slpits or $f$ must be a constant and hence the homotopy class
is trivial.

The above theorem include the following results of \cite{X.Wang, X.Wang1}.

\begin{cor}
Let $(M,g)$ as in theorem \ref{thm-top1}, then $H_{n}(M,{\Z})=0$ and in
particular the conformal infinity is connected.
\end{cor}

As we mentioned in the introduction, all the above application and study are
motivated by the results of \cite{Witten-Yau}. And we are now ready to prove
the promised generalization (theorem \ref{thm-top}) which is an immediate
consequence of theorem \ref{thm-top1}

\noindent\textbf{Proof of theorem \ref{thm-top}.} By Lee's theorem \cite{Lee}%
, $M$ has eigenvalue $\lambda=n^{2}/4$ which is strictly greater than $n-1$
if $n\geq3$. Therefore, by theorem \ref{thm-top1}, the claim follows. If $n=2$, theorem (\ref{thm-top1}) implies that either the claim is true or the manifold splits. As the boundary is assumed to have positive Yamabe invariant, the warped product metric cannot be negative Einstien. Hence the theorem is also true for this case.

\begin{cor}
Let $(M^{n+1},g)$, $n\geq2$, be a conformally compact Einstein manifold with 
$\mbox{Ric}_{g}=-ng$. Suppose that the conformal infinity of $M$ has
positive Yamabe invariant, then $H_{n}(M,{}\mathbb{Z})=0$ and in particular
the conformal infinity is connected.
\end{cor}

\noindent\emph{Proof.} This follows immediately from theorem \ref{thm-top}
and basic results in algebraic topology.

Finally, we apply theorem \ref{thm-top} to prove the a nonexistence theorem
for Einstein manifold in the situation of AdS/CFT correspondence.

\begin{thm}
\label{thm-nonex} Suppose that $N^{n+1}$, $n\ge2$, is a compact manifold
which support a nonpositively curved metric. Then for any embedded disc $%
D\subset N$, the manifold $M=N\setminus D$ has no asymptotically hyperbolic
conformally compact metric of order $C^{1}$ satisfying $\mbox{Ric}\geq-n$
and $\lambda\ge n-1$.
\end{thm}

\noindent\emph{Proof.} It is obvious that $[(M,\partial M), (N, *)]$ is
nontrivial for any point $*\in D$ and $M$ does not split. Therefore, the existence of such metric
contradicts theorem \ref{thm-top}.

\begin{thm}
Suppose that $N^{n+1}$, $n\geq2$, is a compact manifold which support a
nonpositively curved metric. Then for any embedded disc $D\subset N$, the
manifold $M=N\setminus D$ has no conformally compact Einstein metric of
order $C^{3,\alpha}$ satisfying $\mbox{Ric}= -n$.
\end{thm}

\noindent\emph{Proof.} Since $\partial M=\partial D$ has positive Yamabe
invariant, if there exists a conformally compact Einstein metric of order $%
C^{3,\alpha}$ with $\mbox{Ric}=-n$, then Lee's theorem \cite{Lee} implies
that $\lambda=n^{2}/4$ which is greater than $n-1$ for $n\ge3$. This is
impossible by theorem \ref{thm-nonex}. If $n=2$, then the manifold splits with the given form of warped product metric. The negativity of this Einstein metric contradicts the fact that the $\partial M$ has positive Yamabe invariant.

\end{document}